\documentclass[11pt]{article}
\usepackage{amsmath}
\usepackage{amscd}
\usepackage{amsfonts}
\usepackage{amssymb}
\usepackage{amsthm}
\usepackage[OT2,T1]{fontenc}
\DeclareSymbolFont{cyrletters}{OT2}{wncyr}{m}{n}
\DeclareMathSymbol{\Sha}{\mathord}{cyrletters}{"58}

\newcommand{\rar}{\rightarrow}
\newcommand{\calg}{\mathcal}

\newcommand{\strc}{\mathcal{O}_{X}}

\newcommand{\enn}{\text{End}}

\newcommand{\dif}{{\mathcal {D}}\text{iff}}
\newcommand{\cale}{\mathcal{E}}
\newcommand{\difdt}{{{\mathcal {D}}\text{iff}}^{\bullet}}

\newcommand{\diff}{\text{Diff}}
\newcommand{\diffdt}{{\text{Diff}}^{\bullet}}
\newcommand{\hh}[2]{\text{H}^{#1}(X,#2)}
\newcommand{\hoch}[1]{\text{HH}_{#1}}
\newcommand{\choch}[1]{\widehat{\text{HH}}_{#1}}

\newcommand{\ccccyc}[1]{\widehat{\text{HC}}_{#1}}
\newcommand{\compl}{\mathbb C}
\newcommand{\dol}[2]{{\text{K}^{\bullet}}_{{#1 \text{ } #2}}}

\newcommand{\hcoh}[2]{{\mathbb H}^{#1}(X, #2)}
\newcommand{\hooc}[1]{\widetilde{\text{Hoch}}(#1)}
\newcommand{\ccycl}[1]{\widetilde{\text{Cycl}}(#1)}

\newtheorem{thm}{Theorem}
\newtheorem{prop}{Proposition}
\newtheorem{cor}{Corollary}
\newtheorem{lem}{Lemma}
\newtheorem*{thmm}{Theorem A}

\title{Integration over complex manifolds via Hochschild homology.}
\author{Ajay C. Ramadoss}

\begin{document}

\maketitle

\begin{abstract}
Given a holomorphic vector bundle $\cale$ on a connected
compact complex manifold $X$, [FLS]  construct a $\compl$-linear functional $I_{\cale}$ on
$\hh{2n}{\compl}$. This is done by constructing a linear functional on
the $0$-th completed Hochschild homology $\choch{0}{(\dif(\cale))}$ of the sheaf of holomorphic
differential operators on $\cale$ using topological quantum mechanics.
They show that this functional is $\int_X$ if
$\cale$ has non zero Euler characteristic.
They conjecture that this functional is $\int_X$ for all $\cale$. \\

A subsequent work [Ram] by the author proved that the linear
functional $I_{\cale}$ is independent of the vector bundle $\cale$.
This note builds upon the work in [Ram] to prove
that $I_{\cale}=\int_X$ for an arbitrary holomorphic vector bundle
$\cale$ on an arbitrary connected compact complex manifold
$X$. This is done using an argument that is very natural from the geometric point of view. Moreover, this argument enables one to make the
approach to this conjecture developed first in [FLS] and subsequently in [Ram] independent of the Riemann-Roch-Hirzebruch theorem.
This
argument also enables us to extend the construction in [FLS] to a construction of a linear functional
$I_{\cale}$ on $\text{H}^{2n}_{c}(Y,\compl)$ for a holomorphic vector bundle $\cale$ with bounded geometry on an arbitrary connected
 complex manifold $Y$ with bounded geometry and prove that $I_{\cale} = \int_Y$.
  We also generalize a result of [Ram] pertaining to "cyclic homology analogs" of $I_{\cale}$.     \\



\end{abstract}

\section*{Introduction}

Let $X$ be a smooth, connected compact complex manifold. Let $\cale$
be a holomorphic vector bundle on $X$. In this note, the term vector bundle shall
refer to a holomorphic vector bundle over a complex manifold unless explicitly
stated otherwise. Let $\dif(\cale)$ be the
sheaf of holomorphic differential operators on $\cale$. We have a
notion of completed Hochschild homology
$\choch{\bullet}{(\dif(\cale))}$ such that there is an isomorphism
$$\beta_{\cale}: \choch{-i}{(\dif(\cale))} \simeq \hh{2n-i}{\compl}$$ for
every integer $i$. A paper by B.Feigin, A.Losev and B.Shoikhet
[FLS] describes the construction of a $\compl$-linear functional
$\text{tr}$ on $\choch{0}{(\dif(\cale))}$ using topological quantum mechanics.
 Denote the linear functional
$\text{tr} \circ \beta_{\cale}^{-1}$ on $\hh{2n}{\compl}$ by $I_{\cale}$. We call $I_{\cale}$ the
FLS-functional of $\cale$. We also have a notion of completed cyclic homology $\ccccyc{\bullet}{(\dif(\cale))}$
such that $$\ccccyc{-j}{(\dif(\cale))} \simeq \hh{2n-j}{\compl} \oplus \hh{2n-j+2}{\compl} \oplus ... \text{ . }$$
The construction of $\text{tr}$ can be generalized to yield linear functionals $\text{tr}_{2i}$ on $\ccccyc{-2i}{(\dif(\cale))}$
for each $i \geq 0$. The linear functional $\text{tr}_{2i}$ therefore yields a linear functional $I_{\cale,2i,2k}$ on $\hh{2n-2k}{\compl}$ whenever
$0 \leq k \leq i$. \\

[FLS] showed that if $\cale$ is a vector bundle of non-zero Euler
characteristic, then $I_{\cale} = \int_X$ as linear functionals on
$\hh{2n}{\compl}$. This was done using the Riemann-Roch-Hirzebruch
theorem. [FLS] conjectured that $I_{\cale} = \int_X$ in general. We
refer to this conjecture as the integral conjecture for the rest of
this note. [Ram] showed that $I_{\cale} = I_{\calg F}$
 for any two holomorphic vector bundles $\cale$ and $\calg F$ on
 $X$. We may therefore, refer to $I_{\strc}$ as the {\it FLS-functional on }$X$.
 The fact that $I_{\cale}$ is independent of $\cale$ implies
  that $I_{\cale} = \int_X$ as long as $\cale$ is a
 holomorphic vector bundle on a compact complex manifold $X$
 that admits {\it at least one} vector bundle of nonzero Euler
 characteristic. This proved the integral conjecture for compact complex
 manifolds arising out of complex algebraic varieties, since any
 smooth complex algebraic variety has at least one vector bundle of
 nonzero Euler characteristic (see the introduction in [Ram] for an argument proving this assertion).
The integral conjecture for arbitrary compact complex
manifolds however, remained an open question as it is not known
whether or not there exist compact complex manifolds with
 no holomorphic vector bundle of nonzero Euler characteristic. \\

This note proves the integral conjecture in general. This is done
by building upon the work in [Ram] which in turn, was a further
development of the approach to this problem in [FLS]. The argument
used here is very natural from the geometric point of view. Further, it makes the approach to the integral conjecture developed
in [FLS] and subsequently in [Ram] independent of the Riemann-Roch-Hirzebruch theorem.

\textbf{Convention:} Throughout this paper, a connected complex manifold shall mean a connected complex manifold that is complete and has bounded
geometry (positive radius of injectivity + all covariant derivatives of the Ricci curvature are bounded) as a Riemannian manifold. \\

Our approach to the integral conjecture in this note also enables us to extend the
construction of the FLS-functional to the construction of a $\compl$-linear functional $I_{\cale}$ on
$\text{H}^{2n}_{c}(Y,\compl)$ given a vector bundle with bounded geometry (see Section 3.3) $\cale$ on an {\it arbitrary} connected complex manifold
 $Y$ , and prove that $I_{\cale} = \int_Y$ (Theorem 1). Here, $\text{H}^{\bullet}_{c}$ denotes cohomology with compact supports.
 One can also extend the construction of $I_{\cale,2i,2k}$ for vector bundles on compact complex manifolds to vector bundles on
 arbitrary complex manifolds. Given a vector bundle $\cale$ with bounded geometry on an {\it arbitrary} connected complex manifold
 $Y$ with, one can extend the construction of $I_{\cale,2i,2k}$ from the compact complex case to construct a $\compl$-linear functional
  $I_{\cale,2i,2k}$ on $\text{H}^{2n-2k}_{c}(Y,\compl)$. Strengthening a result of [Ram], we show that $I_{\cale,2i,0} = \int_Y $ and $I_{\cale,2i,2k} =0 $
  whenever $k>0$ (Theorem 2).\\

\textbf{Outline of this note.} Section 1 contains certain remarks about the idea used
 in this note. These may help the reader understand the motivation behind this note better.
 Section 2 recalls the construction of the FLS-functional on $X$.
Section 3 proves the integral conjecture in general and extends the construction of the FLS-functional to
vector bundles with bounded geometry on {\it arbitrary} connected complex manifolds. \\

\textbf{Acknowledgements.} I am very grateful to Prof. Boris Tsygan for going through this paper carefully and for his very useful comments and
suggestions. I also thank the referee for going through this paper carefully and for his useful comments and suggestions.

\section{Some remarks.}

\textbf{Remark 1:}We remark that the argument used to prove the integral conjecture in Section 3.1 shows
that to prove the integral conjecture for an arbitrary compact
complex manifold of complex dimension $n$, it suffices to prove it
for {\it one} compact complex manifold of the same complex
dimension. This observation enables us to free this
approach from the Riemann-Roch-Hirzebruch theorem (see Section 3.2). Recall
that the proof of the integral conjecture for vector bundles of
non-zero Euler characteristic used the fact that the class of the
global differential operator $\text{id}$
in $\text{H}^{2n}(X,\compl)$ is $(\text{ch}(\cale).\text{td}_X)_{2n}$ (see [NT1] and [NT2]) along with the Riemann-Roch-\\Hirzebruch theorem.
Theorem 2 of [Ram] together with the argument in Section 3.1 proves that we need this only for {\it one} particular
vector bundle on {\it one} particular compact complex manifold of dimension $n$ for the integral conjecture to hold for
every compact complex manifold of dimension $n$. This is exploited in Section 3.2 to do away with the need for the Riemann-Roch-Hirzebruch
theorem altogether. A special case of the fact that the class of the
global differential operator $\text{id}$
in $\text{H}^{2n}(X,\compl)$ is $(\text{ch}(\cale).\text{td}_X)_{2n}$ is however still used. This fact together with the integral conjecture
implies the Riemann-Roch-Hirzebruch theorem itself, giving yet another proof of the Riemann-Roch-Hirzebruch theorem.\\

\textbf{Remark 2.} The argument in Section 3.1 is also very natural
from the geometric point of view. Let $U$ be an open disc with inclusions into two compact complex manifolds $X$ and $Y$. If $\omega$ is a top
degree differential form on $U$ supported compactly in $U$ then
$$\int_X \omega = \int_U \omega = \int_Y \omega \text{ . }$$ The
crux of this note is to prove "directly" that the Feigin-Losev-Shoikhet linear
functional mimics the above behavior of the integral. This is
exploited along with the results from [Ram] to prove the integral
conjecture in general. We also note that the same idea is behind
the extension of the Feigin-Losev-Shoikhet construction of the integral via topological
quantum mechanics to non-compact complex manifolds as well (Section 3.3). Of
course, cohomology with compact supports has to be used instead of cohomology itself. \\

\textbf{Remark 3.}A related conjecture in [FLS] that has since been
proven in [EnFe] pertained to traces of global holomorphic
differential operators on $\cale$. If $D$ is a global holomorphic
differential operator on $\cale$, $D$ induces endomorphisms on
$\text{H}^i(X,\cale)$ for all $i$. The {\it supertrace} of $D$ ,
$\text{str}(D)$ is given by the formula $$\text{str}(D) = \sum_i
{(-1)}^i \text{tr}(D|_{\text{H}^i(X,\cale)}) \text{ . }$$ Further,
$D$ is seen to yield a class $[D]$ in $\choch{0}{(\dif(\cale))}$
(see [FLS],[Ram]). It follows from the construction of $\text{tr}:\choch{0}{(\dif(\cale))} \rar \compl$
(see [FLS],[Ram]) that $$\text{tr}([D]) = \text{str}(D) \text{ . }$$
Denote the element $\beta_{\cale}([D])$ of $\hh{2n}{\compl}$ by
$[D]$ itself. It was conjectured by [FLS] and proven by [EnFe] that
$$\text{str}([D]) = \int_X [D] \text{ . }$$ We refer to this result as
the supertrace theorem in this note. This result is somewhat similar to Cor 5.6 of [S-S].
Note that the integral conjecture implies the
supertrace theorem. Also note that the supertrace theorem together with Theorem 2 of [Ram]
implies the integral conjecture for any compact complex manifold that admits at least one
holomorphic vector bundle admitting at least one global holomorphic
differential operator with non-zero supertrace. Unfortunately, we do
not know whether every compact complex manifold has this
property. We also point out that by proving the integral conjecture in full generality, this note completes a different "Riemann-Roch
-Hirzebruch theorem free" approach to the
supertrace theorem from that in [EnFe]. Proposition 4.1 of [EnFe] inspired the author to use a "heat kernel"
approach to push the idea outlined in Remark 2 through.\\

\section{Preliminary material}

This section is meant to briefly recall the salient aspects of
earlier work in [FLS] and [Ram]. For further details, the reader may
refer to [FLS] and [Ram]. Let $\text{Dolb}(X,\strc)$ denote the
Dolbeaux resolution of $\strc$. Denote the complex $\dif(\cale)
\otimes_{\strc}\text{Dolb}(X,\strc)$ by $\dif^{\bullet}(\cale)$. Let
$\diffdt(\cale)$ denote the differential graded algebra of global
sections of $\difdt(\cale)$. Let $\dol{}{\cale}$ denote the complex
$\Gamma(X,\cale \otimes_{\strc}\text{Dolb}(X,\strc))$. By basic
Hodge theory (see [Vois] theorem 5.24), $\dol{}{\cale}$ splits into a direct sum of a complex
$\dol{0}{\cale}$ of $\compl$-vector spaces with $0$ differential and
an acyclic complex $\dol{1}{\cale}$. \\

\subsection{The key construction of [FLS]}

The reader may refer ro [FLS] for further details regarding
any assertion made in this section.
 The key construction of [FLS] is of an $A_{\infty}$-morphism $\calg
F$ from $\diffdt(\cale)$ to $\enn(\dol{0}{\cale})$. The
$A_{\infty}$-morphism $\calg F$ induces a map ${\calg
F}_{\text{hoch}}$ from the Hochschild chain complex of
$\diffdt(\cale)$ to that of $\enn(\dol{0}{\cale})$. One thus obtains
a map ${\calg F}_{\text{hoch}*}$ from the Hochschild homology of
$\diffdt(\cale)$ to that of $\enn(\dol{0}{\cale})$. Let
$\hoch{i}(A)$ denote the $i$th Hochschild homology of a graded
algebra $A$. Then,
$$\hoch{i}(\enn(\dol{0}{\cale})) \simeq 0 \text{    } \forall i \neq
0 $$ $$ \hoch{0}(\enn(\dol{0}{\cale})) \simeq \compl \text{ . }$$
The only Hochschild $0$-cycles that have nontrivial images in
$\hoch{0}(\enn(\dol{0}{\cale})) $ are those arising out of degree
$0$ elements of $\enn(\dol{0}{\cale})$. The image in $\compl$ of the
class in $\hoch{0}(\enn(\dol{0}{\cale}))$ of a Hochschild $0$ cycle
arising out of a degree $0$ element $M$ of $\enn(\dol{0}{\cale})$ is
the supertrace $\text{str}(M)$ of $M$. We therefore denote the
identification of $\hoch{0}(\enn(\dol{0}{\cale})) $ with $\compl$ by
$\text{str}$. It follows from this and from the formula (see [FLS])
for ${\calg F}_{\text{hoch}}$ that if $a$ is a degree $k-1$ element
of $\diffdt(\cale)^{\otimes k}$ yielding a Hochschild $0$-cycle of
$\diffdt(\cale)$, then,
\begin{equation} \label{fhoch} \text{tr}({\calg F}_{\text{hoch}*}(a)) =
\sum_{j=0}^{j=k-1} \text{str}({\calg F}_k(\tau^j(a))) \text{ . }
\end{equation} In the above equation, $\tau$ is the $\compl$-endomorphism of
$\diffdt(\cale)^{\otimes k}$ arising out of a cyclic permutation of
factors with the appropriate sign . ${\calg F}_k$ is the $k$ th {\it
Taylor component} of the
$A_{\infty}$-morphism $\calg F$ (for more details, see [FLS]). \\

We now describe the construction of the Taylor components ${\calg
F}_k$ of $\calg F$. \\

\subsubsection{The Taylor components of $\calg F$}

Let $C_k$ denote the configuration space $\{t_1<...<t_k | t_i \in
\mathbb R\}/G^{(1)}$ where $G^{(1)}$ is the one dimensional group of
shifts $(t_1,..,t_k) \rar (t_1+c,...,t_k+c)$. This is a smooth $k-1$
dimensional manifold that is not compact if $k>1$. Note that setting
$\tau_i := t_{i+1} -t_i$ identifies the $C_k$ with the
open orthant $\Pi_{i=1}^{i=k-1} \{\tau_i > 0\} $. Let
$\overline{\{\tau_i > 0\}}$ denote the compactification of $\{\tau_i
\geq 0\}$ by a point at infinity . Let
$\overline{C_k} = \Pi_{i=1}^{i=k-1} \overline{\{\tau_i > 0\}}$. This is a compactification of $C_k$. \\

If $\phi$ is an element of $\enn(\dol{ }{\cale})$, let $[\phi]_i$
denote the endomorphism $id \otimes ... \otimes \phi \otimes ...
\otimes id $ of $\dol{ }{\cale}^{\otimes k}$ where $\phi$ acts on
the $i$th factor from the right. Recall that
$\bar{\partial}_{\cale}^*$ denotes the Hodge adjoint of
$\bar{\partial}_{\cale}$. Similarly, $\Delta_{\cale}$ denotes the
Laplacian of $\bar{\partial}_{\cale}$.

Let $\Phi$ denote the differential form $$[id]_k \circ
[\text{exp}[-d\tau_{k-1}\bar{\partial}_{\cale}^* -
\tau_{k-1}\Delta_{\cale}]]_{k-1} \circ ..... \circ
[\text{exp}[-d\tau_1 \bar{\partial}_{\cale}^* - \tau_1
\Delta_{\cale}]]_{1}
$$ on $C_k$ with values in $\enn(\dol{ }{\cale})^{\otimes k}$ (note
that $\enn(\dol{ }{\cale}^{\otimes k}) \simeq \enn(\dol{
}{\cale})^{\otimes k}$). This extends to a differential form on
$\overline{C_k}$ with values in $\enn(\dol{ }{\cale})^{\otimes k}$.
In addition, there is a composition map from $\enn(\dol{
}{\cale})^{\otimes k}$ to $\enn(\dol{ }{\cale})$ which we shall
denote by $\text{m}_k$. An element $D$ of $\diffdt(\cale)^{\otimes k}$
yields an element of $\enn(\dol{ }{\cale})^{\otimes k}$ which shall
also be denoted by $D$. If ${\calg I}$ and $\Pi$ denote the
inclusion of $\dol{0}{\cale}$ as a direct summand of $\dol{}{\cale}$
and the projection from $\dol{ }{\cale}$ to $\dol{0}{\cale}$
respectively, then
\begin{equation} \label{fk} {\calg F}_k(D) = \int_{\overline{C_k}} \Pi \circ
\text{m}_k(\Phi \circ D) \circ {\calg I} = \Pi \circ
(\int_{\overline{C_k}} \text{m}_k(\Phi \circ D)) \circ {\calg I}
\text{ . }
\end{equation}

That the ${\calg F}_k$ form
the Taylor components of an $A_{\infty}$-morphism was proven in
[FLS]. \\

\subsection{A linear functional on $\choch{0}{(\dif(\cale))}$}

For an open subset $U$ of $X$, let $\diff(\cale)(U)$ and
$\diffdt(\cale)(U)$ denote $\Gamma(U,\dif(\cale))$ and
$\Gamma(U,\difdt(\cale))$ respectively. Let $\text{C}^{\bullet}(\diff(\cale)(U))$ denote
the complex of Hochschild chains of $\diff(\cale)(U)$ (converted into a cochain complex).  We note that the Hochschild
differential on $\text{C}^{\bullet}(\diff(\cale)(U))$ extends to a
differential of degree $1$ on the graded vector space \\ $\oplus_{k \geq 1}
\diff(\cale^{\boxtimes k})(U^k)[k-1]$ where $\cale^{\boxtimes k}$ is
the $k$-fold external tensor power of $\cale$ on $X^k$. We denote
the resulting complex by
$\widehat{\text{C}^{\bullet}(\diff(\cale)(U))}$ . Similarly, we note
that the Hochschild differential on
$\text{C}^{\bullet}(\diffdt(\cale)(U))$ extends to a differential of degree $1$ on
the graded vector space $\oplus_{k \geq 1} \diffdt(\cale^{\boxtimes
k})(U^k)[k-1]$. We denote the resulting complex by
$\widehat{\text{C}^{\bullet}(\diffdt(\cale)(U))}$ . Let
$\widehat{\text{C}^{\bullet}(\dif(\cale))}$ denote the sheaf of
complexes associated to the presheaf $$U \leadsto
\widehat{\text{C}^{\bullet}(\diff(\cale)(U))}$$ of complexes of
$\compl$-vector spaces on $X$. Similarly, let $\hooc{\dif(\cale)}$
denote the sheaf of complexes associated to the presheaf $$U
\leadsto \widehat{\text{C}^{\bullet}(\diffdt(\cale)(U))}$$ of
complexes of
$\compl$-vector spaces on $X$. \\

By definition, $\choch{i}{(\dif(\cale))} =
\hcoh{i}{\widehat{\text{C}^{\bullet}(\dif(\cale))}}$. \\

Unfortunately, $\calg F$ does not automatically yield a map of
complexes from $\hooc{\dif(\cale)}$ to
$\text{C}^{\bullet}({\enn(\dol{0}{\cale})})$ . One however, has the
following facts (see Proposition 6 of [Ram]). Recall
that any $0$-cocycle $\alpha$ of $\text{C}^{\bullet}(\diffdt(\cale))$ is of the form $\sum_k \alpha_k$
where $\alpha_k \in \diffdt(\cale)^{\otimes k}[k-1]$. Note that $\alpha_1 \in \text{Diff}^0(\cale)$ and
$\alpha_k =0$ for almost all $k$. Let $\Pi_0$ denote the projection from $\dol{}{\cale}$ onto the kernel of the Laplacian
$\Delta_\cale:\dol{}{\cale} \rar \dol{}{\cale}$. This is an integral operator with smooth kernel (see [BGV] Chapter 2).  \\

{\it Fact 1:} The linear functionals $$\alpha \mapsto \sum_k \sum_{j=0}^{j=k-1} \text{Str}({\calg F}_k(\tau^j(\alpha_k)))$$
and $$\alpha \mapsto \text{str}(\Pi_0 \alpha_1 \Pi_0) $$ coincide on the space of $0$-cocycles of $\text{C}^{\bullet}(\diffdt(\cale))$.
Denote this linear functional by $I_{\text{FLS}}$. Recall from [FLS] that $I_{\text{FLS}}$ vanishes on $0$-coboundaries. \\

{\it Fact 2:} $I_{\text{FLS}}$ extends to yield a linear functional on the $0$th cohomology of
$\Gamma(X,\hooc{\dif(\cale)})$. We will denote this linear functional by $\widehat{\text{tr}}$. \\

On the other hand, the natural degree preserving map of complexes
from $\widehat{\text{C}^{\bullet}(\dif(\cale))}$ to
$\hooc{\dif(\cale)}$ is a quasiisomorphism since $\diffdt(\cale^{\boxtimes
k})(U^k)$ is quasiisomorphic to $\diff(\cale^{\boxtimes k})(U^k)$ for any $k \geq 1$ and any open
subset $U$ of $X$. Further,
$\hooc{\dif(\cale)}$ is a complex of sheaves of $\compl$-vector
spaces that are modules over the sheaf of smooth functions on $X$.
It follows that the $i$th cohomology of the complex
$\Gamma(X,\hooc{\dif(\cale)})$ is $\hcoh{i}{\hooc{\dif(\cale)}} =
\hcoh{i}{\widehat{\text{C}^{\bullet}(\dif(\cale))}} =
\choch{i}{(\dif(\cale))}$. \\

It follows that $\widehat{\text{tr}}$ is a linear functional on $\choch{0}{(\dif(\cale))}$. Also recall (for instance, [Ram] Lemma 3) that
$\widehat{\text{C}^{\bullet}(\dif(\cale))}$ is quasiisomorphic to the shifted constant sheaf
$\underline{\compl}[2n]$. It follows that $\choch{-i}{(\dif(\cale))} \simeq \hh{2n-i}{\compl}$.
It follows that $\widehat{\text{tr}}$ yields a linear functional on $\hh{2n}{\compl}$ which we denote by $I_{\cale}$.
\\

\section{Generalizing the integral conjecture.}

Let $X$ be a compact complex manifold admitting at least one holomorphic vector
bundle of non-zero Euler characteristic. Let $\cale$ be a
homomorphic vector bundle on $X$. Let $I_{\cale}:\hh{2n}{\compl}
\rar \compl$ be as in the introduction. Since we have already shown
in [Ram] that $I_{\cale}=I_{\calg F}$ for any vector bundle $\calg
F$ on $X$, we may assume without loss of generality that $\cale=
\strc$. Let $\dol{X}{}$ denote the Dolbeaux complex of $\strc$. Let
$C_k$ and $\overline{C_k}$ be as in Section 2.1.1.\\

Let $\dif(X)$ denote $\dif(\strc)$. Choose open discs $U \subset W
\subset  X$. Since \\ $\hooc{\dif(X)}$ is a complex of soft sheaves
quasiisomorphic to $\underline{\compl}[2n]$ , the complex
$\Gamma_c(U,\hooc{\dif(X)})$ is quasiisomorphic to
$\text{H}^{2n+\bullet}_{c}(U,\compl)$. Here, $\Gamma_c$ is the
functor "sections with compact support" and $\text{H}^{\bullet}_{c}$
denotes cohomology with compact support. Note that
$\text{H}^{2n}_{c}(U,\compl) \simeq \compl$ and
$\text{H}^{i}_{c}(U,\compl) =0 \text{  } \forall \text{  } i \neq
2n$.

It follows that $$\text{H}^{0}(\Gamma_c(U,\hooc{\dif(X)})) \simeq
\text{H}^{2n}_{c}(U,\compl)  \simeq \compl \text{ . }$$

Let $[\alpha]_U$ denote the class of $\alpha$ in $\text{H}^{2n}_{c}(U,\compl)$ for any $0$-cocycle $\alpha$
of \\ $\Gamma_c(U,\hooc{\dif(X)})$. We state the following obvious fact as a proposition for emphasis.\\

\begin{prop} There is a Hochschild $0$-cycle $\alpha$ of
$\Gamma_c(U,\hooc{\dif(X)})$ such that $[\alpha]_U \neq 0$ in
$\text{H}^{0}(\Gamma_{c}(U,\hooc{\dif(X)})) \simeq
\text{H}^{2n}_{c}(U,\compl)$.
\end{prop}

Let $\alpha$ be as in Proposition 1. Let $\Delta$ denote the Laplacian of
$\bar{\partial}$ on $\dol{X}{}$. Let $\dol{X}{L^2}$ denote the
Hilbert space of square integrable Dolbeaux forms on $X$. This is a ${\mathbb Z}_2$-graded
Hilbert space. Let $D^k$ denote sheaf associated to
the presheaf $U \leadsto \diffdt(U^k)[k-1]$. Let $\alpha_k$ denote the component of $\alpha$ in $\Gamma_c(U,D^k)$. Note that $\alpha_k = 0$ for almost all $k$.
Note that $\alpha_1$ is a compactly supported element of $\text{diff}^0(U)$. Therefore, $\alpha_1$ may also be thought of as an element
of $\text{diff}^0(X)$.\\

We now recall Proposition 2.45 of [BGV] as a lemma.\\

\begin{lem}
For any scalar $t > 0$, the operator $\alpha_1 \text{e}^{-t\Delta}$ makes sense as a
 trace class operator on
 $\dol{X}{L^2 }$.
\end{lem}

Thus, if $$\varphi(\alpha):= \alpha_1 \text{ , }$$ then
$\varphi(\alpha) \text{e}^{-t\Delta}$ makes sense as a trace class operator
on $\dol{X}{L^2}$. Let $\text{str}_X(\theta)$ denote the supertrace
of $\theta$ for any trace class operator $\theta$ on $\dol{X}{L^2}$. Note that $\alpha$ may
also be thought of as a $0$-cocycle of $\Gamma(X,\hooc{\dif(X)})$. Let $[\alpha]_X$ denote the
class of $\alpha$ in $\text{H}^{2n}(X,\compl)$. If $j_{X}$ denotes the inclusion from $U$ into $X$, then
$[\alpha]_X = j_{X*}[\alpha]_U$. \\

\begin{prop}
$$\lim_{t \rar \infty} \text{str}_X(\varphi(\alpha) \text{e}^{-t \Delta}) =
\int_X [\alpha]_X $$
\end{prop}
\begin{proof}
Let $\dol{0}{X}$ denote the kernel of $\Delta$ .Recall that
$\dol{0}{X}$ is finite dimensional. Recall that $\Delta$ is an
operator on $\dol{X}{L^2}$ with discrete non-negative spectrum that preserves the ${\mathbb Z}_2$-grading. One
can thus find a graded Hilbert space basis of $\dol{X}{L^2}$
 made up entirely of eigenvectors of $\Delta$. Let $\{e_1,....,\}$
 be such a basis with $\lambda_i $ denoting the eigenvalue of $e_i$.
 Let $\langle \text{  ,  }\rangle$ denote the inner product of
 $\dol{X}{L^2}$. We will denote $\varphi(\alpha)$ by $\varphi$ for the rest of this proof.
 \\

 Then, $$\text{str}_X(\varphi \text{e}^{-t \Delta}) = \sum_i \pm
 \langle \varphi \text{e}^{-t\Delta}(e_i),e_i \rangle = \sum_i
 \pm \text{e}^{-{(t-1)}\lambda_i} \langle \varphi \text{e}^{-\Delta} (e_i),e_i \rangle $$
 for any $t >1$ . Note that the above sums converge absolutely by the fact that
 $\varphi \text{e}^{-t \Delta}$ is a trace class operator for any $t>0$. It follows that
 $$\lim_{t \rar \infty} \text{str}_X(\varphi \text{e}^{-t \Delta}) =
 \lim_{t \rar \infty} \sum_i
 \pm \text{e}^{-(t-1)\lambda_i} \langle \varphi \text{e}^{-\Delta}(e_i),e_i \rangle $$
 $$ = \sum_{ \{i | \lambda_i =0\}} \pm \langle \varphi(e_i),e_i \rangle \text{ .
 }$$
 The last sum is a finite sum and is equal to
 $\text{str}(\Pi_{\dol{0}{X}} \circ \varphi \circ {\calg
 I}_{\dol{0}{X}})$ where $\Pi_{\dol{0}{X}}$ and ${\calg
 I}_{\dol{0}{X}}$ are the projection from $\dol{X}{ }$ to
 $\dol{0}{X}$ and the inclusion from $\dol{0}{X}$ into $\dol{X}{ }$
  respectively. Since $X$ has at least one vector bundle of non-zero Euler characteristic,
   $$\text{str}(\Pi_{\dol{0}{X}} \circ \varphi \circ {\calg
 I}_{\dol{0}{X}}) =   \int_X [\alpha]_X$$ by [Ram] Theorem 2. This proves the
 desired proposition.\\
\end{proof}

Note that $$\dol{X}{L^2} = \dol{U}{L^2} \oplus \dol{X \setminus
U}{L^2} $$ as graded Hilbert spaces. Let $\Delta_U$ denote the restriction of $\Delta_X$ to $U$. \\

\begin{prop}
$\varphi(\alpha) \text{e}^{-t \Delta_U}$ is an operator with trace on
$\dol{U}{L^2}$ and $$\text{str}_X(\varphi(\alpha) \text{e}^{-t \Delta_X}) =
\text{str}_U(\varphi(\alpha) \text{e}^{-t \Delta_U})$$ for any $t
>0$.
\end{prop}

\begin{proof}
Denote $\varphi(\alpha)$ by $\varphi$ in this proof. Recall from [BGV] that $\text{e}^{-t\Delta_X}$ is an operator
with smooth kernel $p_t$ (called the heat kernel) and that
$$\text{str}_X(\varphi \text{e}^{-t\Delta_X}) = \int_X \text{str}(\varphi p_t(x,x))|dx| $$
$$= \int_U \text{str}(\varphi p_t(x,x))|dx| = \text{str}_U(\varphi \text{e}^{-t\Delta_U}) \text{ . }$$
The last equality is because the heat kernel on $U$ is unique (see [Don]).
The construction of the heat kernel for a noncompact Riemannian manifold is done by modifying the construction in [BGV] for the compact case.
This modification, as in [Don] goes through, provided the manifold in question has bounded geometry.
The second equality above is because $\varphi$ is a differential operator supported compactly in $U$.\\

\end{proof}
Recall that $[\alpha]_U$ denotes the class of $\alpha$ in
$\text{H}^{2n}_{c}(U,\compl)$.\\

\begin{cor}
$$\lim_{t \rar \infty} \text{str}_U(\varphi(\alpha) \text{e}^{-t \Delta}) =
\int_U [\alpha]_U  $$
\end{cor}

\begin{proof}
Since $\alpha$ is compactly supported on a subset of $U$, $\int_X
[\alpha]_X = \int_U [\alpha]_U$. The corollary now follows from
Proposition 3 and Proposition 2.\\
\end{proof}

\subsection{Proof of the integral conjecture in general.}

Let $Y$ be an arbitrary compact complex manifold with
$\text{dim}_{\compl} Y = \text{dim}_{\compl} X$. Note that we can
find an open disc on $Y$ that we can identify (holomorphically) with $W$. Let $W_X$
and $W_Y$ denote $W$ thought of as open subdiscs of $X$ and $Y$
respectively. Moreover, the Hermitian metric on
$\Omega^{0,\bullet}(W_Y)$ can be chosen such that it coincides with
that on $\Omega^{0,\bullet}(U_X)$ on $U_Y$.

Let $j_X$ and $j_Y$ denote the inclusions from $U$ into $X$ and $Y$
respectively. Let $\alpha$ be as in Proposition 1. Denote $\varphi(\alpha)$ by $\varphi$ in this subsection. Then,
\begin{equation} \label{point1} \int_U [\alpha]_U = \int_X {j_X}_*[\alpha]_U = \int_Y
{j_Y}_*[\alpha]_U \text{ . }\end{equation} On the other hand,
$$\text{str}_X(\varphi \text{e}^{-t \Delta_X}) =
\text{str}_U(\varphi \text{e}^{-t \Delta_U}) = \text{str}_Y(\varphi
\text{e}^{-t \Delta_Y}) $$ for any $t >0$ by Proposition 3. Taking
the limit as $t \rar \infty$ and applying Corollary 1, we get
$$ \lim_{t \rar \infty} \text{str}_Y(\varphi
\text{e}^{-t \Delta_Y}) = \lim_{t \rar \infty} \text{str}_U(\varphi
\text{e}^{-t \Delta_U}) = \int_U [\alpha]_U   \text{ . }$$ By
$\eqref{point1}$,
$$\lim_{t \rar \infty} \text{str}_Y(\varphi
\text{e}^{-t \Delta_Y}) = \int_Y {j_Y}_*[\alpha]_U  \text{ . }$$ But, $[\alpha]_Y = {j_Y}_*[\alpha]_U$
where $\alpha$ on the left hand side is viewed as a $0$-cocycle of
$\Gamma(Y,\hooc{\dif(Y)})$. On the other hand, following the proof of Proposition 2, we see that
$$\lim_{t \rar
\infty} \text{str}_Y(\varphi \text{e}^{-t \Delta_Y}) =
\Pi_{\dol{0}{Y}} \circ \varphi \circ {\calg I}_{\dol{0}{Y}} \text{ .
}$$ The right hand-side is precisely the Feigin-Losev-Shoikhet
linear functional on $Y$ applied to
$[\alpha]_Y$. It follows that the Feigin-Losev-Shoikhet linear
functional on $Y$ applied to $[\alpha]_Y$ is precisely $\int_Y
[\alpha]_Y$. Since $[\alpha]_Y = {j_Y}_*[\alpha]_U \neq 0$, this proves
that the Feigin-Losev-Shoikhet linear
functional on $Y$ is precisely $\int_Y$. \\

\subsection{Proving the integral conjecture without the Riemann-Roch-Hirzebuch theorem.}

Let $Z$ and $Y$ be two compact complex manifolds of complex dimension $n$.The proof in Section 3.1 also proves the following theorem.

\begin{thmm}
The integral conjecture holds for $Z$ iff it holds for $Y$.
\end{thmm}

Consider the vector bundle ${\calg O}_{{\mathbb P}^1}$ on ${\mathbb P}^1_{\compl}$. Consider the differential operator $\text{id}$ on
${\calg O}_{{\mathbb P}^1}$. The following special case of the supertrace theorem uses a result from [NT1] and a hands-on calculation.\\

\begin{prop}
$$1= \text{str}(\text{id}) = \int_{{\mathbb P}^1} [\text{id}] \text{ . }$$
\end{prop}

\begin{proof} Let $T_{{\mathbb P}^1}$ denote the tangent bundle of ${\mathbb P}^1$. That $\text{str}(\text{id})=1$ follows from the fact that $\text{H}^0({\mathbb P}^1,\calg O_{{\mathbb P}^1}) \simeq \compl$
and $\text{H}^0({\mathbb P}^1,\calg O_{{\mathbb P}^1})=0$. By Theorem 7.1.1 of [NT1], $[\text{id}] = (\text{td}(T_{{\mathbb P}^1}))_{2}$.
We therefore need to verify that $\int_{{\mathbb P}^1} (\text{td}(T_{{\mathbb P}^1}))_2 =1$. Note that $T_{{\mathbb P}^1}$ is a line bundle.
It follows that $(\text{td}(T_{{\mathbb P}^1}))_2 = \frac{1}{2}\text{c}_1(T_{{\mathbb P}^1})$ where $\text{c}_1(\cale)$ denotes the first Chern class of
$\cale$. Also, $T_{{\mathbb P}^1} = {\calg O}(2)$. It therefore suffices to show that $\int_{{\mathbb P}^1} \text{c}_1({\calg O}(1)) =1$.
Then, $\text{c}_1({\calg O}(1))$ is the class of the Chern form of ${\calg O}(1)$, which we will denote by $\omega_{\text{ch}}$. Let $z$ denote the local
holomorphic coordinate on an affine line $U \subset {\mathbb P}^1$. Then, $\int_{{\mathbb P}^1} \omega_{\text{ch}} = \int_U \omega_{\text{ch}}$. On the
other hand, on $U$,
$$\omega_{\text{ch}} = \frac{i}{2\pi} \frac{dz \wedge d{\bar{z}}}{(1+|z|^2)^2} $$
by Lemma 3.16 of [Vois]. Setting $z=x+iy$, it follows that
$$\int_U \omega_{\text{ch}} = \frac{1}{\pi} \int_{\mathbb R^2}\frac{dxdy}{(1+x^2+y^2)^2} =1 \text{ . }$$
This proves the desired proposition.

\end{proof}

It follows that the operator $\text{id}^{\otimes n}$ on ${\calg O}_{{{\mathbb P}^1 }^{\times n}}$ also has supertrace $1$.
Further, by Proposition 5 below, after identifying $\text{H}^{2n}({{\mathbb P}^1}^{\times n},\compl)$ with $\text{H}^2({\mathbb P}^1,\compl)^{\otimes n}$,
$$ 1= \int_{{{\mathbb P}^1}^{\times n}} [\text{id}^{\otimes n}] = \int_{{{\mathbb P}^1}^{\times n}} [\text{id}]^{\otimes n} \text{ . }$$
The integral conjecture therefore holds for ${{\mathbb P}^1 }^{\times n}$, and hence (by Theorem A) for any compact complex manifold provided we prove the following
proposition. In the following proposition, $Y$ and $Z$ are compact complex manifolds. $\calg D_1$ and $\calg D_2$ are global holomorphic operators on
$Y$ and $Z$ respectively. As a result ${\calg D}_1 \otimes {\calg D}_2$ is a global holomorphic differential operator on $Y \times Z$.

\begin{prop}
$$ [{\calg D}_1 \otimes {\calg D}_2] = [{\calg D_1}] \otimes [{\calg D}_2]  \text{ . }$$
\end{prop}

\begin{proof}\text{ }\\
 {\it Step 0: Fixing basic notation.}\\
Let $n$ and $m$ denote the complex dimensions of $Y$ and $Z$ respectively. Let $\widehat{\text{C}^{\bullet}(\dif(M))}$ denote the completed Hochschild
chain complex $\widehat{\text{C}^{\bullet}(\dif({\calg O}_M))}$ for any complex manifold $M$ (converted into a cochain complex). We recall from [Bryl] that  $\widehat{\text{C}^{\bullet}(\dif(M))}$
 is quasiisomorphic to the shifted constant sheaf $\underline{\compl}[2d]$  on $M$ where $d$ is the complex dimension of $M$. Denote this quasiisomorphism by
 $i_M$. \\

 Let $\text{C}^{\bullet}(A)$ denote the Hochschild chain complex of an $\compl$-algebra $A$ viewed as a cochain complex.
 Recall that if $A$ and $B$ are $\compl$-algebras , the shuffle product $\Sha$ yields a map of complexes from
  $\text{C}^{\bullet}(A) \otimes \text{C}^{\bullet}(B)$
 to $\text{C}^{\bullet}(A \otimes B)$ (see [Loday], Section 4.2). In particular if $U$ and $V$ are open discs in $Y$ and $Z$ respectively,
 the shuffle product yields a map of complexes from
 $\text{C}^{\bullet}(\text{Diff(U)}) \otimes \text{C}^{\bullet}(\text{Diff}(V))$ to $\text{C}^{\bullet}(\text{Diff}(U \times V))$.
 This further extends to a map of complexes of sheaves of $\compl$-vector spaces on $Y \times Z$ from
 $\widehat{\text{C}^{\bullet}(\dif(Y))} \otimes \widehat{\text{C}^{\bullet}(\dif(Z))}$ to
 $\widehat{\text{C}^{\bullet}(\dif(Y \times Z))}$ which we will denote by $\text{m}_{Sh}$.\\

{\it Step 1: Reduction to a "local check".}\\
 Think of ${\calg D}_1$ and ${\calg D}_2$ as elements of $\Gamma(Y,\widehat{\text{C}^{0}(\dif(Y))})$ and
 $\Gamma(Z,\widehat{\text{C}^{0}(\dif(Z))})$ respectively. Then,
 $$\text{m}_{\text{Sh}}({\calg D}_1 \otimes {\calg D}_2) = {\calg D_1} \otimes {\calg D_2} \in
 \Gamma(Y \times Z,\widehat{\text{C}^{0}(\dif(Y \times Z))}) \text{ . }$$

 To prove this proposition, it therefore suffices to show that the following diagram commutes in the derived category
 $\text{D}_{\text{Sh}_{\compl}}(Y \times Z)$
 of sheaves of $\compl$-vector spaces on $Y \times Z$.\\

 $$\begin{CD}
    \widehat{\text{C}^{\bullet}(\dif(Y))} \otimes \widehat{\text{C}^{\bullet}(\dif(Z))} @> \text{m}_{\text{Sh}} >> \widehat{\text{C}^{\bullet}(\dif(Y \times Z))}\\
    @VV{i_{Y} \otimes i_Z}V      @V{i_{Y \times Z}}VV\\
    \underline{\compl}[2n] \otimes \underline{\compl}[2m] @>>> \underline{\compl}[2n+2m] \\
    \end{CD} $$

Since a sheaf of $\compl$-vector spaces is injective iff it is flasque (see [Riet], Lemma 3.3), the constant sheaf $\underline{\compl}$
is an injective object in the category of $\compl$-vector spaces on $X$. It follows from this that the diagram above commutes in
 $\text{D}_{\text{Sh}_{\compl}}(Y \times Z)$ upto a scalar. Checking that that scalar factor is one is "done locally". Let $U$ and $V$ be open
  discs in $Y$ and $Z$ respectively. It suffices to show that the following diagram commutes upto cohomology in the category
   of complexes of $\compl$-vector spaces.\\

 $$\begin{CD}
    \widehat{\text{C}^{\bullet}(\diff(U))} \otimes \widehat{\text{C}^{\bullet}(\diff(V))} @> \text{m}_{\text{Sh}} >> \widehat{\text{C}^{\bullet}(\diff(U \times V))}\\
    @VV{i_{Y}|_U \otimes i_Z|_V}V      @V{i_{Y \times Z}|_{U \times V}}VV\\
    {\compl}[2n] \otimes {\compl}[2m] @>>> {\compl}[2n+2m] \\
    \end{CD} $$

  The bottom row of the above diagram is the natural identification of $\compl \otimes \compl$ with $\compl$ that takes $1 \otimes 1$ to $1$. \\

{\it Step 2:The "local check".}\\
  Let $\text{sgn}(\sigma)$ denote the sign of a permutation $\sigma \in S_{k}$. Note that if $W$ is any $\compl$-vector space, then
  $\sigma$ acts on $W^{\otimes k}$ on the right as
  $$\sigma(w_1 \otimes ... \otimes w_k) = w_{\sigma(1)} \otimes .... \otimes w_{\sigma(k)} \text{ . }$$

   Let $\omega(z)$ denote the Hochschild 2-cycle
   $$1 \otimes \frac{\partial}{\partial z} \otimes z - 1 \otimes z \otimes  \frac{\partial}{\partial z}+ 1 \otimes 1 \otimes 1$$
   of the Weyl algebra generated by $z$ and $\frac{\partial}{\partial z}$ .
   Then, $$\omega_{2n} := \omega(z_1) \Sha \omega(z_2) \Sha ... \Sha \omega(z_n) $$
   is a Hochschild $2n$-cycle of the Weyl algebra ${\mathcal A}_n$ generated by $z_1,..,z_n$ and
   $\frac{\partial}{\partial z_1},.....,\frac{\partial}{\partial z_n}$.
   If $z_1,...,z_n$ are local holomorphic
  coordinates on $U$, then ${\mathcal A}_n$ is a subalgebra
  of $\text{Diff}(U)$. It follows that $\omega_{2n}$ is a Hochschild $2n$-cycle in $\text{C}^{\bullet}(\text{Diff}(U))$. Note that
  that image of $\omega_{2n}$ in the normalized Hochschild chain complex of $\text{Diff}(U)$ is the normalized Hochschild $2n$-cycle
  $$\sum_{\sigma \in S_{2n}} \text{sgn}(\sigma) 1 \otimes \sigma(\frac{\partial}{\partial z_1}
   \otimes z_1 \otimes ....
  \otimes \frac{\partial}{\partial z_{n}} \otimes z_n)\text{ .} $$
 We recall from [BrGe] and [FT] that
  $$i_Y|_{U} ([\omega_{2n}]) = 1 \text{ . }$$  To check that the above diagram commutes, we only need to verify that
  $$\text{m}_{\text{Sh}} (\omega_{2n} \otimes \omega_{2m}) = \omega_{2n+2m} \text{ . }$$ This is immediate from our definition of
  $\omega_{2n}$.\\

{\it Remark:} \\
The construction of $i_Y|_{U}$ from [Bryl] was what we used in [Ram]. Even with this construction if $i_Y|_{U}$, we can directly verify that
$i_Y|_U([\omega_{2n}])=1$.We now sketch how this can be done. Recall that the cohomology of $\widehat{\text{C}^{\bullet}(\text{Diff}(U))}$ was computed using the
spectral sequence arising out of the filtration induced by a specific filtration $F^{\bullet}$ on $\text{Diff}(U)$. Here, $F^{-k}\text{Diff(U)}$ was the
space of differential operators on $U$ of order at most $k$. Let $z_1,...,z_n$ be local holomorphic coordinates on the cotangent bundle $T^{*}U$
of $U$. Let $y_1,...,y_n$ be local holomorphic coordinates on the fibre of $T^{*}U$. Setting the weight of the $dz_i$ to be $0$ and that of the $
dy_i$ to be $1$ enables us to define the notion of the weight of a holomorphic form on $T^{*}U$. In the next paragraph, differential forms on $T^{*}(U)$
shall always refer to holomorphic differential forms on $T^{*}U$ that are algebraic along the fibres. \\

The $E^{p,q}_{1}$ term of the spectral sequence computing the cohomology of \\$\widehat{\text{C}^{\bullet}(\text{Diff}(U))}$ is precisely the space of
$-p-q$-forms on $T^{*}U$ of weight $-p$ that are algebraic along the fibres. In fact, the image of the cycle $\omega_{2n}$ in $E_1^{-n,-n}$ can be
verified to be the differential form $dy_1 \wedge dz_1 \wedge dy_2 \wedge dz_2 \wedge.... \wedge dy_n \wedge dz_n$. Recall from [Bryl] (Theorem 3.1.1) that the
differential on the $E_1^{\bullet,\bullet}$ terms is the differential of the canonical complex of the Poisson manifold $T^{*}U$.
Moreover the canonical complex of $T^{*}U$ may be identified with the (shifted) De-Rham
complex of $T^{*}U$. Under this identification, $dy_1 \wedge dz_1 \wedge dy_2 \wedge dz_2 \wedge.... \wedge dy_n \wedge dz_n$ is identified with $1$.
This is a De-Rham $0$-cocycle.
It follows that the image of $\omega_{2n}$ in $E_2^{-n,-n}$ is $1$. Since, $E_{2}^{p,q} = 0$ whenever $(p,q) \neq (-n,-n)$, the cohomology of
$\widehat{\text{C}^{\bullet}(\text{Diff}(U))}$ is identified with $E_2^{-n,-n} \simeq \compl$. This verifies that $i_Y|_U([\omega_{2n}])=1$.\\

 \end{proof}

\subsection{The Feigin-Losev-Shoikhet construction for vector bundles on noncompact complex manifolds.}
As a byproduct of this proof, we have in fact, extended the construction of the
Feigin-Losev-Shoikhet linear functionals associated with certain holomorphic bundles on complex manifolds to complex manifolds that are not compact.
Let $\cale$ be a holomorphic vector bundle on an arbitrary connected complex manifold $Y$.
 Let $\Delta_{\cale}$ denote the Laplacian of
$\cale$. This depends on a choice of Hermitian metric for $Y$ as well as for $\cale$. Recall that the Laplacian $\Delta_{\cale} = \Delta^{\cale} + F$ where $\Delta^{\cale}$ is the Laplacian of a connection on $\cale$ (see Definition 2.4 of
[BGV]) and $F \in \Gamma(Y,\text{End}(\cale))$.\\

\textbf{Definition:}We say that $\cale$ has {\it bounded geometry} if for some choice of Hermitian metric on $\cale$, there exists a connection $\triangledown_{\cale}$ on $\cale$ such that
$\Delta_{\cale} = \Delta^{\cale} +F$ where $\Delta^{\cale}$ is the Laplacian of $\triangledown_{\cale}$ and $F \in \Gamma(Y,\text{End}(\cale))$, and
all covariant derivatives of the curvature of $\Delta^{\cale}$ as well as of $F$ are bounded on $Y$.\\

Let $\cale$ be a vector bundle having bounded geometry. Let $\dol{\cale}{L^2}$ denote the (${\mathbb Z}_2$-graded) Hilbert space of square integrable sections of $\dol{\cale}{ }$.
Then, $\text{e}^{-t\Delta_{\cale}}$ can be constructed as an integral operator on $\dol{\cale}{L^2}$
following [Don]. Let $\alpha \in \Gamma_c(Y,\hooc{\dif(\cale)})$.
Let $D^k(\cale)$ denote the sheaf associated with the presheaf $U \leadsto \diffdt(\cale^{\boxtimes k})(U^k)$.
Let $\alpha_k$ denote the component of $\alpha$ in $\Gamma_c(Y,D^k(\cale))$.
Note that $\alpha_k =0$ for almost all $k$.
Let $$\varphi(\alpha)= \alpha_1
\text{ . }$$  The following
proposition generalizes Proposition 6
 of [Ram]. \\

\begin{prop} \text{    }\\
1. $\varphi(\alpha)\text{e}^{-t\Delta_{\cale}}$ makes sense as a trace class operator on $\dol{\cale}{L^2}$ for any $t >0$. \\
2. Further, the map $$\alpha \leadsto \lim_{t \rar \infty} \text{str} (\varphi(\alpha)\text{e}^{-t\Delta_{\cale}})$$ induces a $\compl$-
linear functional on $\text{H}^{0}(\Gamma_c(Y,\hooc{\dif(\cale)}))$. \\
\end{prop}

\begin{proof}
{\it Step 1:}To prove Part 1, note that the support of $\alpha$ can be covered by finitely many open discs $U_1,..,U_m \subset Y$ such that
each $U_i$ is contained in an open disc $W_i \subset Y$ such that $\cale$ is trivial on each $W_i$. Let $Z= Y \setminus \cup_i \overline{U_i}$.
One can find a
partition of unity $\{f_1,...,f_m,f\}$ on $Y$ subordinate to the cover $Y= \cup_i U_i \cup Z$ with $f_i$ supported on $U_i$.
Note that the support of $f_i$ is compact. Also, $f=0$ on the support of $\alpha$.
 Writing $\alpha$ as $\sum_i f_i \alpha$ it suffices to prove part 1 for $\alpha$ compactly supported on an open disc $U$
contained in an open disc $W$ on which $\cale$ is trivial. One may find a compact complex manifold $X$ containing an open
disc $W_X$ with which $W$ can be identified holomorphically. Let $\cale'= {\mathcal O}_{X}^{p}$ where $p$ is the rank of $\cale$.
 The metrics on $X$ and $\cale'$ may be chosen to coincide with those of $U$ and $\cale |_U$ respectively
 on $U_X$. Let $\alpha_X$ denote $\alpha$ thought of as an element of $\Gamma(X, \hooc{\dif(\cale')})$.
 By an easy generalization of Lemma 1,  $\varphi(\alpha_X)\text{e}^{-t\Delta_{\cale'}}$ makes sense as a trace class
 operator on $\dol{\cale'}{L^2}$ for any $t >0$. Since $\alpha_X$ is supported on a compact subset of $U_X$ and since
 $\dol{\cale'}{L^2} = \dol{\cale'|_{U_X}}{L^2} \oplus \dol{\cale'|_{X \setminus U_X}}{L^2}$, an easy extension of Proposition 3 implies
  that   $\varphi(\alpha_X)\text{e}^{-t\Delta_{\cale'}}$ makes sense as a trace class
 operator on $\dol{\cale'|_{U_X}}{L^2}$ for any $t >0$ and that $$\text{str}_X(\varphi(\alpha_X)\text{e}^{-t\Delta_{\cale'}})
 = \text{str}_{U_X}(\varphi(\alpha_X)\text{e}^{-t\Delta_{\cale'}}) \text{ . }$$ But since $\alpha_X$,$U_X$ and $\Delta_{\cale'} |_{U_X}$
 are identified with $\alpha$,$U$ and $\Delta_{\cale} |_U$ respectively, $\varphi(\alpha)\text{e}^{-t\Delta_{\cale}}$ makes sense as a
 trace class operator on $\dol{\cale|_U}{L^2}$ for any $t >0$ and
 $$\text{str}_U( \varphi(\alpha)\text{e}^{-t\Delta_{\cale}}) = \text{str}_X(\varphi(\alpha_X)\text{e}^{-t\Delta_{\cale'}})\text{ . }$$
 for any $t>0$. Noting that $\dol{\cale}{L^2} = \dol{\cale|_U}{L^2} \oplus \dol{\cale|_{Y \setminus U}}{L^2}$ and noting that
 $\alpha$ is supported on a compact subset of $U$, we see that  $\varphi(\alpha)\text{e}^{-t\Delta_{\cale}}$ makes sense as a
 trace class operator on $\dol{\cale}{L^2}$ for any $t >0$ and
 \begin{equation} \label{genprop3} \text{str}_Y( \varphi(\alpha)\text{e}^{-t\Delta_{\cale}}) =
 \text{str}_X(\varphi(\alpha_X)\text{e}^{-t\Delta_{\cale'}}) \end{equation}
 for any $t>0$. This proves part 1.\\

{\it Step 2:}  By $\eqref{genprop3}$,
$$\lim_{t \rar \infty} \text{str}_Y (\varphi(\alpha)\text{e}^{-t\Delta_{\cale}}) =
\lim_{t \rar \infty} \text{str}_X (\varphi(\alpha_X)\text{e}^{-t\Delta_{\cale'}}) \text{ . }$$ Moreover, since $X$ is compact, a trivial modification
of the argument proving Proposition 2 will show that the right hand side is finite . It follows that
$\alpha \leadsto \lim_{t \rar \infty} \text{str} (\varphi(\alpha)\text{e}^{-t\Delta_{\cale}})$ yields a
linear functional on the space of compactly supported sections of degree $0$ of $\hooc{\dif(\cale)}$ that are supported on any fixed open disc
$U \subset W$ such that $W \subset Y$ is an open disc on which $\cale$ is trivial.
 That this
extends to a linear functional on the space of degree $0$ elements of $\Gamma_c(Y,\hooc{\dif(\cale)})$ follows from a partition of unity argument similar
to that used to prove part 1 of this proposition.\\

{\it Step 3:} To show that $\alpha \leadsto \lim_{t \rar \infty} \text{str} (\varphi(\alpha)\text{e}^{-t\Delta_{\cale}})$ yields a linear functional on
$\text{H}^{0}(\Gamma_c(Y,\hooc{\dif(\cale)}))$ , we need to show that\\
$\lim_{t \rar \infty} \text{str} (\varphi(d_{\text{hoch}}\beta)\text{e}^{-t\Delta_{\cale}}) = 0$ for any degree $-1$ element $\beta$ of\\
$\Gamma_c(Y,\hooc{\dif(\cale)})$. Here, $d_{\text{hoch}}$ is the differential of the complex\\ $\Gamma_c(Y,\hooc{\dif(\cale)})$. Once more,
as in the proof of Part 1, one can first show that it suffices to show that
$\lim_{t \rar \infty} \text{str} (\varphi(d_{\text{hoch}}\beta)\text{e}^{-t\Delta_{\cale}}) = 0$ for any degree $-1$ element $\beta$ of
$\Gamma_c(Y,\hooc{\dif(\cale)})$ supported on a subset of an open disc $U \subset W$ such that $W \subset Y$ is an open disc
on which $\cale$ is trivial. In this case, if $X$ and $\cale'$
are as in Step 1 of this proof, then $$\lim_{t \rar \infty} \text{str}_Y (\varphi(d_{\text{hoch}}\beta)\text{e}^{-t\Delta_{\cale}})
 = \lim_{t \rar \infty} \text{str}_X (\varphi(d_{\text{hoch}}\beta_X)\text{e}^{-t\Delta_{\cale'}}) \text{ . } $$
 Since $X$ is compact, by an easy generalization of proposition 2,
 $$ \lim_{t \rar \infty} \text{str}_X (\varphi(d_{\text{hoch}}\beta_X)\text{e}^{-t\Delta_{\cale'}}) =
 \Pi_{\dol{0}{\cale'}} \circ  \varphi(d_{\text{hoch}}\beta_X) \circ {\calg I}_{\dol{0}{\cale'}} \text{ . }$$
 The right hand side is precisely $ \int_X [d_{\text{hoch}} \beta_X] = 0 $ by the integral conjecture for compact
 complex manifolds. This proves part 2 of the desired proposition.\\

\end{proof}

Note that $\hooc{\dif(\cale)}$ is a complex of soft sheaves that are modules over the sheaf of smooth functions on $Y$. Further, $\hooc{\dif(\cale)}$
 is quasiisomorphic to $\widehat{\text{C}^{\bullet}(\dif(\cale))}$ which in turn is quasiisomorphic to the shifted constant sheaf $\underline{\compl}[2n]$
 (see [Ram] , Lemma 3). It follows that \\$\text{H}^{0}(\Gamma_c(Y,\hooc{\dif(\cale)})) \simeq \text{H}^{2n}_{c}(Y,\compl)$. By Part 2 of Proposition 6,
 we have constructed a $\compl$-linear functional on $\text{H}^{2n}_{c}(Y,\compl)$, which we will denote by $I_{\cale}$.
 The formula for this linear functional on $\text{H}^{2n}_{c}(Y,\compl)$
  coincides with that for the FLS functional on $\cale$ as constructed in [Ram] when $Y$ is compact. The following generalization of Theorem 2 of [Ram]
   holds.

   \begin{thm} Let $\cale$ be a holomorphic vector bundle having bounded geometry on a connected complex manifold $Y$. Then,
   $$I_{\cale} = \int_Y :\text{H}^{2n}_{c}(Y,\compl) \rar \compl \text{ . } $$
   \end{thm}

   \begin{proof}
   Let $U \subset W \subset Y$ be open discs. Choose a $0$-cocycle $\alpha$ of \\ $\Gamma_c(Y,\hooc{\dif(\cale)})$ such that $[\alpha]_U \neq 0$
   in $\text{H}^{2n}_{c}(U,\compl)$. Let $U_X$, $X$ and $\cale'$ be as in Step 1 of the proof of Proposition 6.
   Then, $$I_{\cale'}(j_{X*}[\alpha]_U) = \int_{X} j_{X*}[\alpha]_U $$ by the integral conjecture for compact complex manifolds. But,
    $\int_X j_{X*}[\alpha]_U =
   \int_U [\alpha]_U  = \int_{Y} j_{Y*}[\alpha]_U $. Moreover,
   $$ I_{\cale}(j_{Y*}[\alpha]_U) = \lim_{t \rar \infty} \text{str}_Y(\varphi(\alpha)\text{e}^{-t\Delta_{\cale}}) $$ $$
   I_{\cale'}(j_{X*}[\alpha]_U) =  \lim_{t \rar \infty} \text{str}_X(\varphi(\alpha_X)\text{e}^{-t\Delta_{\cale'}}) \text{ . }$$
   By $\eqref{genprop3}$ , $I_{\cale}(j_{Y*}[\alpha]_U) = I_{\cale'}(j_{X*}[\alpha]_U)$.
    It follows that $\int_{Y} j_{Y*}[\alpha]_U = I_{\cale}(j_{Y*}[\alpha]_U)$. Since
    $j_{Y*}[\alpha]_U \neq 0$ and $\text{H}^{2n}_{c}(Y,\compl)$ is a one dimensional $\compl$-vector space, the desired theorem follows.

   \end{proof}

\textbf{Remark:} One may define the "completed Hochschild homology with compact support" of $\dif(\cale)$. The completed Hochschild homology
$\widehat{\text{HH}}^{c}_{\bullet}(\dif(\cale))$ of $\dif(\cale)$ is given by $\text{R}\Gamma_c(\widehat{\text{C}^{\bullet}(\dif(\cale))})$. Note that
since $\hooc{\dif(\cale)}$ is a complex of soft sheaves on $X$ quasiisomorphic to
$\widehat{\text{C}^{\bullet}(\dif(\cale))}$, $$ \widehat{\text{HH}^{c}_{i}(\dif(\cale))} \simeq \text{H}^i(\Gamma_c(Y,\hooc{\dif(\cale)}))\text{ . }$$
Also note that if $Y$ is compact, $\Gamma_c=\Gamma$. It follows that $ \widehat{\text{HH}^{c}_{i}(\dif(\cale))} = \widehat{\text{HH}_{i}(\dif(\cale))}$
if $Y$ is compact. It now becomes easy to observe that Proposition 6 extends the definition of $I_{\cale}$ to vector bundles on arbitrary manifolds by generalizing
the definition of the linear functional $\text{tr}$ on $\choch{0}{(\dif(\cale))}$ in the compact case to that of a linear functional on
$ \widehat{\text{HH}^{c}_{0}(\dif(\cale))}$. \\

Let $\widehat{\text{CC}^{\bullet}(\diffdt(\cale)(U))}$ denote the completed Tsygan's double complex of $\diffdt(\cale)(U)$. Denote the sheafification
of the presheaf $$U \leadsto \text{tot}(\widehat{\text{CC}^{\bullet}(\diffdt(\cale)(U))})$$ by $\ccycl{\dif(\cale)}$. This is a (soft) sheaf of modules over
the sheaf of smooth functions on $Y$ that is quasiisomorphic to the completed cyclic complex of $\dif(\cale)$. Since the latter complex is
quasiisomorphic to $\underline{\compl}[2n] \oplus \underline{\compl}[2n+2] \oplus ..... $, $$\text{H}^{-2i}(\Gamma_c(Y,\hooc{\dif(\cale)})) \simeq
\text{H}^{2n-2i}_{c}(Y,\compl) \oplus \text{H}^{2n-2i+2}_{c}(Y,\compl) \oplus ... \oplus \text{H}^{2n}_{c}(Y,\compl) \text{ . }$$

One may also note that a $-2i$ cocycle of $\Gamma_c(Y,\ccycl{\dif(\cale)})$ arises out of a tuple $(\beta_{-2i},....,\beta_{0},..,\beta_{l})$ where
$\beta_k \in \Gamma_c(Y,\hooc{\dif(\cale)}^k)$ if $k$ is even and $\beta_k \in \Gamma_c(Y,\widetilde{\text{bar}(\dif(\cale))}^k)$ if $k$ is odd. Note that
the terms of the "bar complex" $\widetilde{\text{bar}(\dif(\cale))}$ are the same as those of $\hooc{\dif(\cale)}$ but the differential of the complex
$\widetilde{\text{bar}(\dif(\cale))}$ is an extension of the bar differential rather than the Hochschild differential. The proof of the following proposition,
which uses Proposition 6 , is very similar to that of Proposition 13 of [Ram] and is thus omitted.\\

\begin{prop}
 The map
 $$(\beta_{-2i},....,\beta_{0},..,\beta_{l}) \leadsto \lim_{t \rar \infty} \text{str} (\varphi(\beta_0)\text{e}^{-t\Delta_{\cale}})$$
 induces a $\compl$-linear functional on $\text{H}^{-2i}(\Gamma_c(Y,\hooc{\dif(\cale)}))$.
 \end{prop}

 We have therefore, constructed a $\compl$- linear functional on
  $\text{H}^{2n-2i}_{c}(Y,\compl) \oplus \text{H}^{2n-2i+2}_{c}(Y,\compl) \oplus ... \oplus \text{H}^{2n}_{c}(Y,\compl)$.
 We will denote the composition of this $\compl$-linear functional with the inclusion of $\text{H}^{2n-2k}_{c}(Y,\compl)$ into
 $\text{H}^{2n-2i}_{c}(Y,\compl) \oplus \text{H}^{2n-2i+2}_{c}(Y,\compl) \oplus ... \oplus \text{H}^{2n}_{c}(Y,\compl) $ as a
 direct summand by $I_{\cale,2i,2k}$ whenever $0 \leq k \leq i$. The following generalization of Theorem 3 of [Ram] holds. Its proof completely
 parallels the proof of Theorem 3 in [Ram] Section 5. \\

 \begin{thm}Let $\cale$ be a holomorphic vector bundle having bounded geometry on a complex manifold $Y$. Then,
 $$I_{\cale,2i,0} = \int_Y : \text{H}^{2n}_{c}(Y,\compl) \rar \compl$$  for any $i \geq 0$. Further,
 $$I_{\cale,2i,2k} = 0 \text{  } \forall \text{  } {0 < k \leq i} \text{ . } $$
 \end{thm}

\textbf{Remark:} One may define the "completed cyclic homology with compact support" $\widehat{\text{HC}^{c}_{\bullet}(\dif(\cale))}$. By definition,
$\widehat{\text{HC}^{c}_{\bullet}(\dif(\cale))} = \text{R}\Gamma_{c}(Y,\widehat{\text{Cycl}(\dif(\cale))})$ where $\widehat{\text{Cycl}(\dif(\cale))}$
is the completed cyclic chain complex of $\dif(\cale)$. Since $\ccycl{\dif(\cale)}$ is quasiisomorphic to $\widehat{\text{Cycl}(\dif(\cale))}$, and is
a complex of soft sheaves on $Y$, $$\widehat{\text{HC}^{c}_{i}(\dif(\cale))} \simeq
\text{H}^{i}(\Gamma_c(Y,\ccycl{\dif(\cale)})) \text{ . } $$ $I_{\cale, 2i,2k}$ is constructed in the non-compact case by observing that the construction
of $\text{tr}_{2i}$ in the compact case generalizes to yield a linear functional on $ \widehat{\text{HC}^{c}_{-2i}(\dif(\cale))}$ for each $i \geq 0$ by
Proposition 7.\\

\end{document}